\documentclass[letterpaper, 10 pt, conference]{ieeeconf}  

\IEEEoverridecommandlockouts                              

\usepackage{pbox}
\usepackage{enumerate}
\usepackage{mathtools}
\usepackage{dsfont}
\usepackage{bm}
\usepackage{bbold}
\def\PR {\mathop{\mathds{P}\kern0pt}}
\def\EXP {\mathop{\mathds{E}\kern0pt}}
\usepackage{algorithmicx}
\usepackage[ruled]{algorithm}
\usepackage{algpseudocode}
\usepackage{ntheorem}
\newtheorem{lemma}{Lemma}
\newtheorem{corollary}{Corollary}
\newtheorem{theorem}{Theorem}
\newtheorem{problem}{Problem}
\newtheorem{remark}{Remark}
\usepackage{graphics}
\usepackage{epsfig} 
\usepackage{amsmath} 
\usepackage{amssymb}  
\usepackage{epstopdf}
\usepackage{mathrsfs}
\usepackage{algorithmicx}
\usepackage[ruled]{algorithm}
\usepackage{algpseudocode}
\newcommand{\Exp}[1]{\mathbb{E}[#1]}
\newcommand{\Prob}[1]{\mathbb{P}(#1)}   

\usepackage{verbatim}
\usepackage{float}
\usepackage{ntheorem}

\newtheorem{definition}{Definition}
 
\newcommand{\argmin}{\operatornamewithlimits{argmin}}

\renewcommand\thesection{\arabic{section}}
\renewcommand\thesubsection{\thesection.\arabic{subsection}}

\title{\LARGE \bf
Reinforcement Learning in Decentralized Stochastic Control Systems  with Partial History Sharing 
}

\author{Jalal Arabneydi$^{1}$ and Aditya Mahajan$^{2}$
\thanks{$^{1}$Jalal Arabneydi is with the Department of Electrical Engineering, 
        McGill University, 3480 University St., Montreal, QC, Canada 
        {\tt\small jalal.arabneydi@mail.mcgill.ca}}%
\thanks{$^{2}$Aditya Mahajan is with the Department of Electrical Engineering, McGill University, 3480 University St., Montreal, QC, Canada
        {\tt\small aditya.mahajan@mcgill.ca}}%
}

\begin{document}

\maketitle

\vspace*{-5.2cm}{\footnotesize{Proceedings of American Control Conference, 2015.}}
\vspace*{4.45cm}
\thispagestyle{empty}
\pagestyle{empty}

\begin{abstract}
In this paper, we are interested in systems with multiple agents that wish to collaborate in order to accomplish a common task while a) agents have different information (decentralized information) and b) agents do not know the  model of the system completely i.e., they may know the model partially or may not know it at all. The agents must learn the optimal strategies by interacting with their environment i.e., by decentralized Reinforcement Learning (RL). The presence of multiple agents with different information makes decentralized reinforcement learning conceptually more difficult than centralized reinforcement learning. In this paper, we develop a decentralized reinforcement learning algorithm that learns $\epsilon$-team-optimal solution for  partial history sharing information structure, which encompasses a large class of decentralized control systems including  delayed sharing, control sharing, mean field sharing, etc.  Our approach consists of two main steps. In the first step, we convert the decentralized control system to an equivalent centralized POMDP (Partially Observable Markov Decision Process) using an existing approach called common information approach. However, the resultant POMDP requires the complete knowledge of system model. To circumvent this requirement, in the second step, we introduce a new concept called ``Incrementally Expanding Representation'' using which we construct a finite-state RL algorithm whose approximation error converges to zero exponentially fast.  We illustrate the proposed approach and verify it numerically  by obtaining a decentralized Q-learning algorithm for two-user Multi Access Broadcast Channel (MABC) which is a benchmark example for decentralized control systems.

\end{abstract}

\section{INTRODUCTION}
\vspace*{-.0cm}
\subsection{Motivation}

Decentralized decision making is relevant in a wide range of applications
ranging from networked control systems, robotics,
transportation networks, communication networks, sensor networks, and economics. There is a 
rich history of research on optimal stochastic control of decentralized system.
We refer the reader to \cite{Yuksel2013stochastic} for a detailed review.

Most of the literature assumes that the system model is completely known to all decision
makers; however, in practice, such knowledge may only  be available partially or may  not be available. Hence, it is crucial for decision makers  to be able to  learn the optimal solutions. In the literature,   learning in centralized stochastic control  is well studied and there exist  many approaches  such as model-predictive
control, adaptive control, and reinforcement learning.  This is in contrast to the learning in decentralized stochastic control; it is not immediately
clear on how centralized learning approaches would work for decentralized systems. In this paper, we propose a novel Reinforcement Learning (RL) algorithm for a class of decentralized stochastic control systems that guarantees team-optimal solution.
%
%

 Existing approaches for multi-agent learning may be categorized as follows: exact methods and heuristics.
The exact methods rely on the assumption that the information structure is such
that all agents can consistently update the Q-function. These include
approaches that rely on social convention and rules to restrict the decisions
made by the agents~\cite{ Spaan:2002}; approaches that
use communication to convey the decisions to all agents~\cite{Vlassis:2003}; and approaches that assume that the Q-function decomposes into a
sum of terms, each of which is independently updated by an
agent~\cite{ Kok:2005}.
Heuristic approaches include joint action learners heuristic~\cite{Claus:2006}, where each
agent learns the empirical model of the system in order to estimate the
control action of other agents; frequency maximum Q-value
heuristic~\cite{Kapetanakis:2002}, where agents keep track of the frequency
with which each action leads to a ``good'' outcome; heuristic
Q-learning~\cite{Matignon:2007}, which assigns a rate of punishment for each
agent; and distributed Q-learning~\cite{Huang:2005}, which uses
predator-prey models to assign heuristic sub-goals to individual agents. 
\nocite{Busoniu2008comprehensive} To the best of our knowledge, there is no RL approach that guarantees
team-optimal solution. In this paper, we present such an approach.

We describe the system model and problem formulation in Sections~\ref{Main_system_model} and~\ref{Main_problem_formulation}, respectively. We state the main challenges  in Section~\ref{Main_difficulties} and  our main contributions in Section~\ref{Contribution}.  In Section~\ref{section:PHS}, we present a brief  preliminary on Partial History Sharing (PHS) information structure. We  describe our  approach in two basic steps in Section~\ref{Approach}. In Section~\ref{Decentralized Implementation}, we mention a few points on the  implementation. Based on the proposed approach, we develop a RL algorithm for a benchmark example  with numerical results in Section~\ref{MABC_example}.

\vspace*{-.1cm}

\subsection{Notation}
We use upper-case letters to denote random variables (e.g. $X$) and lower-case letters to denote their realizations (e.g. $x$). We use the short-hand notation $X_{a:b}$  for the vector $(X_a,X_{a+1},\ldots,X_{b})$ and bold letters to denote vectors e.g. $\mathbf Y=(Y^1,\ldots,Y^n)$.  $\mathbb{P}(\cdot)$ is the probability of an event, $\mathbb{E}[\cdot]$ is the expectation of a random variable, and $|\cdot|$ is the absolute value of a real number. $\mathbb{N}$ refers to the set of natural numbers and $\mathbb{Z}^+=\mathbb{N}\cup \{0\}$.
\vspace*{-.1cm}

\subsection{System Model} \label{Main_system_model}
Let $X_t \in \mathcal{X}$ denote the state of a dynamical system controlled by $n$ agents. At time $t$, agent $i$ observes $Y^i_t \in \mathcal{Y}^i$ and chooses $U^i_t \in \mathcal{U}^i$. For ease of notation, we denote the joint actions and the joint observations by $\mathbf U_t=(U^1_t,\ldots,U^n_t) \in \mathcal{U}$ and $\mathbf Y_t=(Y^1_t,\ldots,Y^n_t) \in \mathcal{Y}$, respectively. The dynamics of the system are given by \vspace{-.2cm}
\begin{equation}
X_{t+1}=f(X_t,\mathbf{U}_t,W^s_t),\vspace{-.2cm}
\end{equation} 
and the observations are given by \vspace{-.2cm}
\begin{equation}
\mathbf{Y}_t=h(X_t,\mathbf{U}_{t-1},W^o_t),\vspace{-.2cm}
\end{equation}
where  $\{W^s_t\}_{t=1}^\infty$ is an i.i.d. process with probability distribution function $P_{W^s}$,  $\{W^o_t\}_{t=1}^\infty$ is an i.i.d. process  with probability distribution function $P_{W^o}$, and $X_1$ is the initial state  with probability distribution function $P_{X}$. The primitive random variables $\{X_1, \{W^s_t\}_{t=1}^\infty, \{W^o_t\}_{t=1}^\infty\}$ are mutually independent and defined on a common probability space.
 
For ease of exposition, we assume all system variables are finite valued. Let $I^i_t \subseteq \{\mathbf Y_{1:t}, \mathbf{U}_{1:t-1}\}$ be information available at agent $i$ at time $t$. The collection $(\{I^i_t\}_{t=1}^\infty, i=1,\ldots,n)$ is called the \textit{information structure}. In this paper, we restrict attention to an information structure called \textit{partial history sharing}  \cite{Nayyar2013CIA}, which will be defined later.

At time $t$, agent $i$ chooses action $U^i_t$ according to \textit{control law} $g^i_t$ as follows\vspace{-.1cm}
\begin{equation}
U^i_t=g^i_t(I^i_t).\vspace{-.1cm}
\end{equation}
We denote $\bm{g}^i=(g^i_1,g^i_2,\ldots)$ as \textit{strategy} of agent $i$ and $\bm g=(\bm g^1,\ldots,\bm g^n)$ as joint strategy of all the agents. The performance of strategy $\mathbf{g}$ is measured by the following infinite-horizon discounted cost
\begin{equation}\label{Total_cost}
J(\bm g)=\mathbb{E}^{\bm g}\left[ \sum_{t=1}^\infty \beta^{t-1} \ell( X_t, \mathbf U_t) \right],
\end{equation}
where $\beta \in (0,1)$ is the discount factor, $\ell$ is the per-step cost function, and the expectation is with respect to a joint probability distribution on $(X_{1:\infty}, \mathbf{U}_{1:\infty})$ induced by the joint probability distribution on the primitive random variables and the choice of strategy $\bm g$.

A strategy $\bm g^\ast$ is optimal if for any other strategy $\bm g$,
$J(\bm{g}^\ast) \leq J(\bm g).$
For  $\epsilon\hspace{-.05cm}>0$,   strategy $\bm g^\ast\hspace{-.05cm}$ is $\epsilon$-optimal, if for any other strategy~$\hspace{-.05cm} \bm g$,
$J(\bm{g}^\ast) \leq J(\bm g)+ \epsilon.$

\vspace{-.0cm}

\subsection{Problem Formulation}\label{Main_problem_formulation}

We will consider three different setups that differ in the assumptions about the knowledge of the model. For all the setups, we will assume that the action and the observation spaces as well as the information structure, the discount factor $\beta$, and an upper-bound on the per-step cost are common knowledge between all agents. The setups differ in the assumptions about state space $\mathcal{X}$, system dynamics and observations ($f,h$), probability distributions ($P_X,P_{W^s},P_{W^o}$), and cost structure $\ell$. These include two setups, 1) complete-knowledge of the  model, and  2) incomplete-knowledge of the model which includes two sub-cases: 2a) partial-knowledge of the model and 2b) no-knowledge of the model.

In general, the complete-knowledge of the model is required to  find an optimal strategy $\bm{g}^\ast$. However, in practice, there are many applications where such information is not completely available or is not available at all. In such applications, the agents must learn the optimal strategy by interacting with their environment. This  is known as reinforcement learning (RL). If the agents have partial knowledge of the model, the setup is called \textit{model-based} RL. If the agents have no knowledge of the model, setup is called \textit{model-free}~RL.

Define $L:=\max_{x,\mathbf u} |\ell(x, \mathbf u)|$. We are interested in the following problem.
\begin{problem}\label{Problem1}
Given the information structure,  action spaces $\{\mathcal{U}^i\}_{i=1}^n$, observation spaces $\{\mathcal{Y}^i\}_{i=1}^n$, discount factor $\beta$, the upper-bound $L$ on  per-step cost, and any $\epsilon>0$,  develop a (model-based or model-free)  reinforcement learning algorithm  using which the agents  learn  an $\epsilon$-optimal strategy~$\bm g^\ast$.
\end{problem}
\vspace*{-.2cm}
\subsection{Main Difficulties}\label{Main_difficulties}

Given the complete knowledge of system model, finding team-optimal solution in decentralized control systems is  conceptually challenging due to the decentralized nature of information available to the agents. The agents need to cooperate with each other to fulfill a common objective while they have different perspectives about themselves, other agents, and the environment. This discrepancy in perspectives makes establishing cooperation among agents difficult; we refer reader to \cite{Murphey2002cooperative2002} for details. Thus, finding team-optimal solution is even more challenging when agents have only partial knowledge or no knowledge of system model.  Hence, it is difficult to  \emph{consistently} learn strategies in such settings.
\vspace*{-.2cm}
%
\subsection{Contributions}\label{Contribution}
Below, we mention our main contributions in this paper.

1) We propose a novel  approach to perform reinforcement learning  in  a large class of  decentralized stochastic control systems  with partial history sharing (PHS) information structure that guarantees $\epsilon$-team-optimal solution. In particular, our approach combines the  \textit{common information approach} of \cite{Nayyar2013CIA} with any RL algorithm of Partially Observable Markov Decision Processes (POMDP). The approach works in two steps. In the first step, the common information approach is used  to convert the decentralized control problem to an equivalent centralized POMDP and in the second step, a RL algorithm is used to provide a learning scheme to identify an $\epsilon$-optimal strategy in the resultant POMDP.  Note that any RL algorithm of POMDPs may be used in the second step; however, we develop a new methodology for the second step as explained below.

2) We propose a novel methodology to perform reinforcement learning in \textit{centralized} POMDPs as an intermediate step of the two-step approach described above. (This methodology by itself may be of interest due to the fact that developing RL algorithm in POMDPs is difficult). The methodology consists of three parts: 1) converting the POMDP to a countable-state MDP $\Delta$ by defining a new concept that we call Incrementally Expanding Representation (IER), 2) approximating $\Delta$ with a sequence of finite-state MDPs $\{\Delta_N\}_{N=1}^\infty$, and 3) using a RL algorithm to learn an optimal strategy of MDP $\Delta_N$.  We show  that the performance of the RL strategy converges to the optimal performance exponentially as  $N \rightarrow \infty$. We use this methodology in the second step of the two-step approach.

3)  Using the proposed two-step approach, we develop a RL algorithm for two-user Multi Access Broadcast Channel (MABC) which is  used as a benchmark for decentralized control systems. Numerical simulations validate that the RL algorithm converges to an optimal strategy.

\section{Preliminaries on Partial History Sharing }\label{section:PHS}
Herein, we  present a simplified version of  partial history sharing information structure, originally presented in~\cite{Nayyar2013CIA}.

\begin{definition}[\cite{Nayyar2013CIA}, Partial History Sharing (PHS)] 
 Consider a decentralized control system with $n$ agents. Let $I^i_t$ denote the information available to agent $i$ at time $t$. Assume $I^i_t \subseteq I^i_{t+1}$. Then,  split the information at each agent into two parts: \textit{common information} $C_t= \bigcap_{i=1}^n I^i_t$ i.e. the information shared between all agents and \textit{local information} $M^i_t=I^i_t \backslash C_t$ that is the local information of agent $i$.  Define $Z_t:=C_{t+1}\backslash C_{t}$ as common observation, then $C_{t+1}=Z_{1:t}$. An information structure is called \textit{partial history sharing} when the following conditions are satisfied:
\begin{enumerate}
\item [a)] The update of local information \[M^i_{t+1} \subseteq \{M^i_t,U^i_t,Y^i_{t+1}\}\backslash Z_t, \quad  i \in \{1,\ldots,n\}.\]

\item [b)] For every agent $i$, the size of the local information~$M^i_t$ and the size of the common observation $Z_t$ are uniformly bounded in time $t$.
\end{enumerate}
\end{definition}

These conditions are fairly mild and are satisfied by a large class of models. Examples include delayed sharing \cite{NayyarDelayed:2011}, periodic sharing \cite{ooi1997separation},   mean-field sharing \cite{JalalMFsharing:2014}, etc. Even for models that do not satisfy the above conditions directly, it is often possible to identify sufficient statistics that satisfy the above conditions, e.g., control sharing \cite{Aditya2013controlsharing}.

\begin{remark}\label{pure-decentralized_information}
Note that the conditions (a) and (b) are valid even if there is no common information between agents i.e., $C_t=\emptyset$. Hence, the decentralized control systems with  pure decentralized information  (i.e. no information commonly shared) falls into PHS information structure.
\end{remark}

\section{Approach}\label{Approach}
In this part, we derive a RL algorithm for systems with PHS information structure.  Our approach consists of  two steps. In the first step, we consider the setup of the complete-knowledge of the model and  use the \textit{common information approach} of \cite{Nayyar2013CIA}  to convert the decentralized control problem to  an equivalent centralized POMDP.  In the second step, we consider the setup of incomplete-knowledge of the model and develop a finite-state RL algorithm based on the POMDP obtained in the first step.

%

  \subsection{Step 1: An Equivalent Centralized POMDP}\label{CIA_general}
  In this section,  we  present common information approach of \cite{Nayyar2013CIA} and its main results for the setup of complete-knowledge of the  model  described in Section \ref{Main_system_model}.

Let $\mathcal{M}^i$ and $\mathcal{Z}$ denote the spaces of realizations of local information of agent $i$ and common observation, respectively. Consider a virtual \textit{coordinator} that observes the common information $C_t$ shared between all agents  and chooses $(\Gamma^1_t,\ldots,\Gamma^n_t)$, where $\Gamma^i_t:\mathcal{M}^i \mapsto \mathcal{U}^i$ is the mapping from the local information of agent $i$ to action of agent $i$ at time~$t$, according to
 \begin{equation}
 \Gamma^i_t= \psi^i_t(C_t), \quad  i \in \{1,\ldots,n\} .  \vspace{-.1cm}
\end{equation} 
We call $\bm \psi_t:=\{\psi^1_t,\ldots,\psi^n_t\}$  the \textit{coordination law} and  $\bm \Gamma_t=(\Gamma^1_t,\ldots,\Gamma^n_t)$ the \textit{prescription}. The agents use this prescription to choose their actions  as follows:
\begin{equation}
 U^i_t=\Gamma^i_t(M^i_t),\quad  i \in \{1,\ldots,n\} . \vspace{-.1cm}
\end{equation}
We denote the space of  mappings $\Gamma^i_t$ by $\mathcal{G}^i$ and the space of prescriptions $\mathbf \Gamma_t$ by $\mathcal{G}=\prod_{i=1}^n \mathcal{G}^i$. In the sequel, for ease of notation, we will use the following compact form for the coordinator's law,\vspace{-.1cm}
\begin{equation}
\bm \Gamma_t=\bm \psi_t(C_t).\vspace{-.1cm}
\end{equation}

We call $\bm \psi=\{\bm \psi_1,\bm \psi_2,\ldots\}$ as the \textit{coordination strategy}. In the \textit{coordinated system}, dynamics and cost function are as same as those in the original problem in Section \ref{Main_system_model}. In particular, the infintie-horizon discounted cost in the coordinated system is as follows:
\begin{equation}
J(\bm \psi)=\mathbb{E}^{\bm \psi}\left[ \sum_{t=1}^\infty \beta^{t-1} \ell(\mathbf X_t, \bm \Gamma^1_t(M^1_t),\ldots,\Gamma^n_t(M^n_t)) \right].
\end{equation}

\begin{lemma}[\cite{Nayyar2013CIA}, Proposition 3]\label{CIA_equivalence}
The original system described in Section \ref{Main_system_model} with PHS information structure is  equivalent to the coordinated system.
\end{lemma}
We denote $\bm M_t=(M^1_t,\ldots,M^n_t)$ as the joint local information. According to \cite{Nayyar2013CIA}, $\Pi_t=\Prob{ X_t, \mathbf M_t |Z_{1:t-1},\bm \Gamma_{1:t-1}}$ is an information state for the coordinated system. It is shown in~\cite{Nayyar2013CIA} that: 

\begin{enumerate}
\item There exists a function $\phi$ such that\vspace{-.1cm}
\begin{equation}\label{POMDP_update_function}
\Pi_{t+1}=\phi(\Pi_t,\bm \Gamma_t, Z_{t}).\vspace{-.1cm}
\end{equation}
\item The observation $Z_t$ only depends on $(\Pi_t,\Gamma_t)$ i.e.\vspace{-.1cm}
\begin{equation}\label{POMDP_observation_transition}
\Prob{Z_t |\Pi_{1:t},\bm \Gamma_{1:t}}=\Prob{Z_t|\Pi_{t},\bm \Gamma_{t}}.\vspace{-.1cm}
\end{equation} 
\item There exists a function $\hat{\ell}$ such that\vspace{-.1cm}
\begin{equation}\label{POMDP_cost}
\hat{\ell}(\pi_t,\hspace{-.07cm}\bm \gamma_t)\hspace{-.1cm}=\hspace{-.1cm}\Exp{\ell(X_t,\hspace{-.05cm}\mathbf U_t|Z_{1:t-1}\hspace{-.1cm}=\hspace{-.1cm}z_{1:t-1}\hspace{-.05cm},\hspace{-.05cm}\bm \Gamma_{1:t}\hspace{-.1cm}=\hspace{-.1cm}\bm \gamma_{1:t})}.\vspace{-.1cm}
\end{equation}
\end{enumerate} 

Assume that the initial state $\pi_1$ is fixed. Let $\mathcal{R}$ denote the reachable set of above centralized POMDP  that contains all the realizations of $\pi_t$  generated by $\pi_{t+1}=\phi(\pi_t,\bm \gamma, z), \forall \bm \gamma \in \mathcal{G}, \forall z \in \mathcal{Z},\forall t\in \mathbb{N}$,  with initial information state $\pi_1$. Note that  since all the variables are finite valued, then $\mathcal{G}$ (set of all prescriptions $\bm \gamma$) and $\mathcal{Z}$ (set of all observations of the coordinator) are finite sets. Hence, $\mathcal{R}$ is at most a countable~set.

\begin{theorem}[\cite{Nayyar2013CIA}, Theorem 5]\label{Theorem1}
Let $\hspace{-.05cm}\bm  \psi^\ast(\pi)$ be any argmin of the right-hand side of  following dynamic program. For $\pi \hspace{-.1cm} \in \hspace{-.1cm} \mathcal{R}$,
\begin{equation*} \label{CIA_dp}
V(\pi)=\min_{\bm \gamma}(\hat{\ell}(\pi,\bm \gamma)+\beta \Exp{V(\phi(\pi,\bm \gamma,Z_{t}))|\Pi_t=\pi,\bm \Gamma_t= \bm \gamma}),\vspace{-.1cm}
\end{equation*}
where $\bm \gamma=(\gamma^1,\ldots,\gamma^n)$ and the minimization is over all functions $\gamma^i \in \mathcal{G}^i, i \in \{1,\ldots,n\}$. Then, the joint  stationary strategy $\bm g^\ast=(\bm g^{1,\ast},\ldots,\bm g^{n,\ast})$ is optimal such that
\begin{equation*}
g^{i,\ast}(\pi, m^i):=\bm \psi^{i,\ast}(\pi)(m^i), \quad  \pi \in \mathcal{R}, m^i \in \mathcal{M}^i, \forall i. 
\end{equation*}

\end{theorem}


In the next step, we develop a finite-state RL algorithm based on the obtained POMDP for the setup of  incomplete-knowledge of the model.

\subsection{Step 2: Finite-State RL Algorithm For POMDP}

In the previous step, we identified a centralized POMDP that is equivalent to the decentralized control system with PHS information structure. However, the obtained POMDP requires  the complete knowledge of the model. To circumvent this requirement, we introduce a new concept that we call \textit{Incrementally Expanding Representation} (IER). The main feature of IER is to remove the dependency of the POMDP from the complete knowledge of the model. Based on a proper IER, in this step, we  develop a finite-state RL algorithm. This step consists of three parts. In part (1), we convert the POMDP to a countable-state MDP $\Delta$  without loss of optimality. In part (2), we construct a sequence of finite-state MDPs $\{\Delta_N\}_{N=1}^\infty$ of MDP $\Delta$. In part (3), we use a generic RL algorithm  to learn an optimal strategy of $\Delta_N$.

\begin{definition}[Incrementally Expanding Representation]\label{IER}
Let $\{\mathcal{S}_{k}\}_{k=1}^\infty$ be a sequence of finite sets such that  $\mathcal{S}_{1} \subsetneq \mathcal{S}_2 \subsetneq \ldots \subsetneq \mathcal{S}_k \subsetneq \ldots $, and $\mathcal{S}_1$ is a singleton, say $\mathcal{S}_1=\{s^\ast\}$. Let $\mathcal{S}=\lim_{k \rightarrow \infty} \mathcal{S}_k$ be the countable union of above finite sets, $B: \mathcal{S} \rightarrow \mathcal{R}$ be a surjective function that maps $\mathcal{S}$ to the reachable set $\mathcal{R}$, and $\tilde{f}: \mathcal{S} \times \mathcal{G} \times \mathcal{Z} \rightarrow \mathcal{S}$.
The tuple $\langle \{\mathcal{S}_{k}\}_{k=1}^\infty, B, \tilde{f}\rangle$ is called an \emph{Incrementally Expanding Representation} (IER), if it satisfies the following properties: 

(P1) \textit{Incremental Expansion}: For any $\bm \gamma \in \mathcal{G}, z \in \mathcal{Z},$ and $s\in \mathcal{S}_k$, we have that \vspace{-.1cm}
\begin{equation}\label{IER_P1}
\tilde{f}(s,\bm \gamma,z) \in \mathcal{S}_{k+1}. \vspace{-.1cm}
\end{equation}

(P2) \textit{Consistency}: For any $(\bm \gamma_{1:t-1},z_{1:t-1})$, let $\pi_t$ and $s_t$  \vspace{-.1cm}
 be the states obtained by recursive application of \eqref{POMDP_update_function} and \eqref{IER_P1} starting from $\pi_1$ and $``s^\ast$, respectively. Then,
 \begin{equation}\label{IER_B}
 \pi_t=B(s_t).  \vspace{-.1cm}
 \end{equation}
\end{definition}

In general, every decentralized control system with PHS information structure has at least one IER. In the following example, we present a generic IER  that is valid for every system with PHS information structure.

\emph{Example 1:} Let $S_1=\{ \emptyset \}$, $S_2=\{\emptyset\} \cup \{\mathcal{G} \times \mathcal{Z}\}$, and   $S_{k+1}=S_k \cup 
\{\mathcal{G} \times \mathcal{Z}\}^k, k \in \mathbb{N}$. Let $S=\lim_{k \rightarrow \infty} S_k$ and $B: \mathcal{S} \rightarrow \mathcal{R}$ such that  \vspace{-.1cm}
\begin{equation*}
B(\emptyset)\hspace{-.05cm}=\hspace{-.05cm}\pi_1,  B(s_{k+1})\hspace{-.05cm}=\hspace{-.05cm}\phi(\phi(....,\bm\gamma_{k-1},z_{k-1}),\bm \gamma_k,z_{k})\hspace{-.05cm}=\hspace{-.05cm}\pi_{k+1}, \vspace{-.1cm}
\end{equation*}
where $s_{k+1}=((\bm \gamma_1,z_1),\ldots,(\bm \gamma_k,z_{k})) \in \mathcal{S}_{k+1}$. Define $\tilde{f}$ as  follows:
\begin{equation*}
\tilde{f}(s,\bm \gamma,z)=s \circ \bm{\gamma} \circ z,
\end{equation*}
where $\circ$ denotes concatenation. By construction,  tuple $\langle \{\mathcal{S}_{k}\}_{k=1}^\infty, B, \tilde{f}\rangle$ satisfies (P1) and (P2), and hence is an IER.

\subsubsection{Countable-state MDP $\Delta$}\label{step2:A}
Let the tuple $\langle \{\mathcal{S}_{k}\}_{k=1}^\infty, B, \tilde{f}\rangle$  be an IER of the POMDP obtained in the first step. Then, define MDP $\Delta$ with countable
state space $\mathcal{S}$, finite action space $\mathcal{G}$, and dynamics $\tilde{f}$ such that:

\textbf{(F1)}  The initial state is singleton $s^\ast$. The state $S_t \in  \mathcal{S}_k$, $k \leq t$, evolves as follows: for  $\bm \Gamma_t \in \mathcal{G}, Z_t \in \mathcal{Z},$
\begin{equation}\label{Delta_dynamics}
S_{t+1}=\tilde{f}(S_t,\bm \Gamma_t,Z_t), \quad  S_{t+1} \in \mathcal{S}_{k+1}
\end{equation}
where observation $Z_t$ only depends on $(S_t,\bm \Gamma_t)$ (that is a consequence of \eqref{POMDP_observation_transition} and consistency property in \eqref{IER_B}).
At time $t$,  there is a cost depending on the current state $S_t \in \mathcal{S}$ and action $ \bm \Gamma_t \in \mathcal{G}$ given by\vspace{-.1cm}
\begin{equation}\label{Countable_state_cost}
\tilde{\ell}(S_t,  \bm \Gamma_t):=\hat{\ell}(B(S_t),  \bm \Gamma_t)=\hat{\ell}(\Pi_t,\bm \Gamma_t).\vspace{-.1cm}
\end{equation}

\textbf{(F2)} State space $\mathcal{S}$, action space $\mathcal{G}$, and dynamics $\tilde{f}$ do not depend on the unknowns.

The performance of a stationary strategy $ \tilde{\bm \psi}: \mathcal{S} \mapsto \mathcal{G}$ is quantified by \vspace{-.2cm}
\begin{equation}\label{Delta_total_cost}
\tilde{J}(\tilde{\bm \psi})=\mathbb{E}^{\tilde{\bm \psi}}\left[\sum_{t=1}^\infty \beta^{t-1} \tilde{\ell}(S_t,\bm \Gamma_t) \right].
\end{equation}

There may exist more than one IER that satisfy above features. For instance, the IER of Example~1 always satisfies (F1) and  (F2) (that is model-free). This IER  can also be used in the model-based cases; however, in the model-based cases, due to having partial knowledge of the model, one may be able to find a simpler IER. See Section \ref{MABC_example} for an example.

\begin{lemma}\label{Delta_lemma}
Let $\tilde{\bm \psi^\ast}$ be an optimal strategy for MDP $\Delta$. Construct a strategy $\bm \psi^\ast$ for the coordinated system as follows:
\begin{equation}
\tilde{\bm \psi^\ast}(s)=:\bm \psi^\ast(B(s)), \quad
 \forall s \in \mathcal{S}.
\end{equation} 
Then, $\tilde{J}(\tilde{\bm \psi^\ast})=J(\bm \psi^\ast)$ and $\bm \psi^\ast$ is an optimal strategy for the coordinated system, and therefore can be used to generate an optimal strategy for the decentralized control system.
\end{lemma}
Proof is omitted due to lack of space.

\subsubsection{Finite-state incrementally expanding MDP $\Delta_N$ }
In this part, we construct a series of finite-state MDPs $\{\Delta_N\}_{N=1}^\infty$, that approximate the countable-state MDP $\Delta$ as follows. Let $\Delta_N$ be a finite-state MDP with state space $\mathcal{S}_N$ and action space $\mathcal{G}$. The transition probability of $\Delta_N$ is constructed as follows. Pick any arbitrary set $D^\ast \in S_N$. Remap every transition in $\Delta$ that takes the state $s \in \mathcal{S}_N$ to $s'\in \mathcal{S}_{N+1}\backslash \mathcal{S}_N$ to a transition from $s \in \mathcal{S}_N$ to any (not necessarily unique) state  in $D^\ast$. In addition, the per-step cost function of $\Delta_N$ is simply a restriction of~$\tilde{\ell}$ to~$\mathcal{S}_N \times \mathcal{G}$.

%

%
%
%

We assume that there exists an action or a sequence of actions that if taken, the system  transmits to a known state $d^\ast$ in~$D^\ast$. For example, suppose there is a reset action  in the system.  After executing the reset action, the state of the system is reset and transmitted to a known state $d^\ast \in D^\ast$.
%
 Let $\tau_N \in \mathbb{N}$ be the longest amount of time during which $S_t$, $t \leq \tau_N$, stays in $\mathcal{S}_N$  under dynamics $\tilde{f}$, optimal strategy $\tilde{\bm \psi^\ast}$, and any arbitrary sample path of $z_{1:\tau_N-1}$, i.e., 
 \begin{equation}
 S_{t}=\tilde{f}(S_{t-1},\tilde{\bm \psi^\ast}(S_{t-1}),Z_{t-1}) \in \mathcal{S}_N, \quad  \forall t \leq \tau_N. 
 \end{equation} 
Let $\tilde{\bm \psi_N^\ast}$ and $\tilde{J}_N(\tilde{\bm \psi_N^\ast})$ be an optimal stationary strategy  of $\Delta_N$ and  the optimal cost (performance) of $\Delta_N$, respectively.

\begin{theorem}\label{Theorem2}
The difference in performance between $\Delta$ and $\Delta_N$ is bounded as follows: 

\begin{equation}
 |\tilde{J}(\tilde{\bm \psi^\ast})-\tilde{J}_N(\tilde{\bm \psi_N^\ast})|\leq \frac{2\beta^{\tau_N}}{1-\beta}L.
\end{equation}

\end{theorem}
Proof is omitted due to lack of space.

The upper-bound provided in Theorem \ref{Theorem2} requires knowledge on $(\tilde{f}$, $\tilde{\bm \psi^\ast}$,$\mathcal{Z})$. However, according to \eqref{Delta_dynamics}, $\tau_N$ is always  equal or greater than $N$ i.e. $N \leq \tau_N$. Hence, one can obtain a more conservative error-bound (larger upper-bound) than the error-bound (upper-bound) in Theorem \ref{Theorem2} that does not require any knowledge on  $(\tilde{f}$, $\tilde{\bm \psi^\ast}$,$\mathcal{Z})$ as follows. 

\begin{corollary}\label{Corollary_conservative_bound}
The difference in performance between $\Delta$ and $\Delta_N$ is bounded as follows: \vspace{-.2cm}
\begin{equation}
 |\tilde{J}(\tilde{\bm \psi^\ast})-\tilde{J}_N(\tilde{\bm \psi_N^\ast})|\leq \frac{2\beta^{N}}{1-\beta}L.\vspace{-.0cm}
\end{equation}
\end{corollary}

\subsubsection{Finite-state RL algorithm} \label{Finite-state RL algorithm}
Let $\mathcal{T}$ be a generic (model-based or model-free) RL algorithm  designed for finite-state MDPs with infinite horizon  discounted  cost. By a generic RL algorithm, we mean any algorithm which fits to the following framework.  At each iteration $k \in \mathbb{N}$, $\mathcal{T}$ knows the state of system, selects one action, and observes an instantaneous cost and the next state.  The strategy learned (generated) by $\mathcal{T}$ converges to an optimal strategy as $k \rightarrow \infty$.

 Let $\mathcal{T}$ operate on MDP $\Delta_N$ such that, at iteration $k$, it knows the state of the system $s_k \in \mathcal{S}_N$, selects one action $\bm \gamma_k \in \mathcal{G}$, and observes an instantaneous cost $\ell_k$ (which is a realization of the incurred cost $\ell(X_k, \mathbf U_k)$ at the original decentralized system).  According to \eqref{POMDP_cost} and \eqref{Countable_state_cost}, we have 
\begin{equation}
\Exp{\ell(X_k,\mathbf U_k)|S_{1:k},\bm \Gamma_{1:k}}=\tilde{\ell}(S_k,\bm \Gamma_k), \quad S_k \in \mathcal{S}_{N}.\vspace{-.1cm}
\end{equation} 
Hence, the instantaneous cost $\ell_k$ may be interpreted as a realization of the per-step cost of $\Delta_N$. Given dynamics  $\tilde{f}$, $\mathcal{T}$ observes $z_k \in \mathcal{Z}$ and computes the next state 
$s_{k+1}=\tilde{f}(s_k,\bm \gamma_k, z_k).$
If $s_{k+1} \in \mathcal{S}_{N+1} \backslash \mathcal{S}_N$, then an action (or a sequence of actions) that transmits the state of system to a known state in $S_N$ i.e. $s_{k+1}=d^\ast \in D^\ast$ will be taken; otherwise, the system will continue from $s_{k+1} \in \mathcal{S}_N$. 

Let $\tilde{\bm \psi}_N^k:\mathcal{S}_N \rightarrow \mathcal{G}$ be the learned strategy associated with RL algorithm $\mathcal{T}$ operating on  MDP $\Delta_N$ at iteration $k$. Then,  $\mathcal{T}$ updates its strategy  $\tilde{\bm \psi}^{k+1}_N$ based on the observed  cost $\ell_k$ and the transmitted next state $s_{k+1}$ by executing action $\bm \gamma_k$ at state $s_k$. We assume $\mathcal{T}$ converges to an optimal strategy $\tilde{\bm \psi^\ast_N}$ as $k \rightarrow \infty$ such that \vspace{-.1cm}
\begin{equation}\label{Generic_RL_convergence}
\lim_{k \rightarrow \infty} |\tilde{J}_N(\tilde{\bm \psi}_N^k) -\tilde{J}_N(\tilde{\bm \psi_N^\ast})|=0. \vspace{-.1cm}
\end{equation} 
Now, we need to convert (translate) the strategies in $\Delta_N$ to strategies in the original decentralized control problem described in Section~\ref{Main_system_model}, where the actual learning happens. Hence, we define a strategy $\bm g_N^k:=(g_N^{k,i},\ldots,g_N^{k,n})$, at iteration $k$, as follows:
\begin{equation}\label{Optimal_learned_strategy}
\hspace{-.06cm} g^{k,i}_N(s,m^i):=\tilde{\bm \psi}_N^{k,i}(s)(m^i), \forall s \in \mathcal{S}_N,\forall m^i \in \mathcal{M}^i,\forall i ,
\end{equation} 
where $\tilde{\bm \psi}_N^{k,i}$ denotes the $i$th term of $\tilde{\bm \psi}_N^{k}$. 


\vspace*{-0cm}
\vspace{-0.0cm}%
\alglanguage{pseudocode}
\begin{algorithm}[b!]
\small
\caption{\hspace{0.35cm}Finite-State RL Algorithm}
\label{alg.1}
\begin{algorithmic}[1]
\State  Given $\epsilon>0$, choose a sufficiently large $N \in \mathbb{N}$ such that $\frac{2\beta^N}{1-\beta}L \leq \epsilon$. Then, construct state space $\mathcal{S}_N$, action space $\mathcal{G}$, and dynamics $\tilde{f}$. Initialize $s_1=s^\ast$. 

\State    At iteration $k \in \mathbb{N}$, RL algorithm $\mathcal{T}$ picks  $\bm \gamma_k=(\gamma^1_k,\ldots,\gamma^n_k) \in \mathcal{G}$ at state $s_k \in \mathcal{S}_N$. 
Then, agent $i\in \{1,\ldots,n\}$ takes action $u^i_k$ according to the chosen prescription $\gamma^i_k$ and local information $m^i_k \in \mathcal{M}^i$ as follows:\vspace{-.1cm}
\begin{equation*}
u^i_k=\gamma^i_k(m_k^i), \quad \forall i.\vspace{-.1cm}
\end{equation*}

\State  Based on the taken actions, the system incurs a cost $\ell_k$, evolves, and   generates  new information i.e. $(\{m^i_{k+1}\}_{i=1}^n,z_k)$. Every agent $i$ observes $z_k$ because it is common observation. Based on $z_k$, all agents  consistently compute the next state\vspace{-.1cm}
\begin{equation*}
 s_{k+1}=\tilde{f}(s_k,\bm \gamma_k,z_k).\vspace{-.1cm}
\end{equation*}
 If $s_{k+1} \notin  \mathcal{S}_N$, then agents take an action (or a sequence of actions) that transmits the state of system to a state $s_{k+1}=d^\ast \in \mathcal{S}_N$; otherwise, the system proceeds from $s_{k+1} \in \mathcal{S}_N$. Note that during the reset process, the algorithm is paused till the system lands in a state in $S_N$.
 
 \State $\mathcal{T}$ updates its strategy from $\tilde{\bm \psi}_N^{k}$  to $\tilde{\bm \psi}_N^{k+1}$ based on performing action $\bm \gamma_k$ at state $s_k$ and transmission  to next state $s_{k+1}$ with instantaneous cost $\ell_k$. 
 
  \State \noindent $k \leftarrow k+1$, and go to step 2 until termination.
\Statex
\end{algorithmic}
  \vspace{-0.2cm}%
\end{algorithm}
\vspace*{-.0cm}%

\begin{theorem}\label{Theorem3}\label{Theorem3}
Let $J^\ast$ be  the optimal performance of the original decentralized control system given in \eqref{Total_cost}. Then, the approximation  error associated with using the learned strategy is bounded as follows:
\begin{equation}
 \lim_{k \rightarrow \infty}|J^\ast-J(\bm g^k_N)|=|\tilde{J}(\tilde{\bm \psi^\ast})-\tilde{J}_N(\tilde{\bm \psi_N^\ast})|\leq \epsilon_N,
\end{equation}
where $\epsilon_N= \frac{2\beta^{\tau_N}}{1-\beta}L\leq  \frac{2\beta^{N}}{1-\beta}L$. Note that  the error goes to zero exponentially in $N$.
\end{theorem}
The proof follows from Theorem \eqref{Theorem1}, Lemma \ref{Delta_lemma}, Theorem \ref{Theorem2}, and Corollary \ref{Corollary_conservative_bound}.

\begin{remark}
In general, a finite-state RL algorithm similar to  Algorithm~\ref{alg.1} can be derived for centralized POMDPs since we do not impose any restriction on the POMDP in step~2. The only required assumption is the existence of an action (or a sequence of actions) that prevents the system to transmit to information states that are generated after a sufficiently long time.  The  existence of such a  reset strategy (``reset button") or an approximate reset strategy (``homing strategy'') is a standard assumption in the literature. See \cite{pivazyan2002polynomial,EyalPhD:2005} and references therein.
\end{remark}


\section{Decentralized Implementation}\label{Decentralized Implementation}

All  agents are provided with state space $S_N$, action space $\mathcal{G}$, and dynamics $\tilde f$ as described in Section~\ref{step2:A}. Note that to obtain above knowledge, every agent must only know the information structure of system,  action spaces~$\{\mathcal{U}^i\}_{i=1}^n$, observation spaces~$\{\mathcal{Y}^i\}_{i=1}^n$, discount factor $\beta$,  upper-bound~$L$ on  per-step cost, and  $\hspace{-.05cm}\epsilon>0$.  In addition, agents may have partial knowledge of the model of system or  may not. 

Agents observe the instantaneous cost of the system. They have access to a common shared random number generator for the purpose of exploring the system consistently.  Given state space $\mathcal{S}_N$, action space $\mathcal{G}$, and dynamics $\tilde{f}$, Algorithm~\ref{alg.1} can be executed in a distributed manner because every agent can independently run Algorithm \ref{alg.1}; agreeing upon a deterministic rule to
break ties while using argmin ensures that all agents
compute the same optimal strategy. Note that no more information needs to be shared; hence, \emph{no communication} is required. According to Remark~\ref{pure-decentralized_information}, Algorithm~\ref{alg.1} also works for the  pure decentralized control systems, when there is no information commonly shared between agents.

Suppose that the generic RL algorithm $\mathcal{T}$ in Section \ref{Finite-state RL algorithm} is Q-learning. Then, in off-line learning, every agent is allowed to have different step sizes (independently from the step sizes of other agents). However, in on-line learning, since we also need to consistently exploit the system, the step sizes should be chosen consistently, e.g., based on the number of visit to  pair of state and  prescription, i.e.,~$(s,\bm \gamma)$.

\section{Example: MABC}\label{MABC_example}

In this section, we provide an example to illustrate our approach. In this example, we consider the setup of partial knowledge of the model.
\subsection{Problem Formulation}

Consider a two-user multiaccess broadcast system. At time $t$, $W^i_t \in
\{0,1\}$ packets arrive at each user according to independent Bernoulli
processes with $\PR(W^i_t = 1) = p^i \in (0,1)$, $i=1,2$. Each user may store only $X^i_t
\in \{0,1\}$ packets in a buffer. If a packet arrives when the user-buffer is full,
the packet is dropped. Both users may transmit $U^i_t \in \{0,1\}$ packets over a shared broadcast
medium. A user can transmit only if it has a packet, thus $U^i_t \le X^i_t$. If
only one user transmits at a time, the transmission is successful and the
transmitted packet is removed from the queue. If both users transmit
simultaneously, packets ``collide'' and remain in the queue. Thus, the state
update for users~1 and 2 is:\vspace{-.2cm}
\begin{equation}
X^i_{t+1} = \min( X^i_t - U^i_t+ U^1_tU^2_t + W^i_t, 1),\quad i=1,2 \vspace{-.2cm}
\end{equation}

Due to the broadcast nature of the communication channel, each user observes the transmission decision of the other user i.e. information at each user at time~$t$ is~ $(X^i_{t},\mathbf{U}_{1:t-1}),~i\in \{1,2\}$. 
Each user chooses a transmission decision as
\begin{equation}
U^i_t = g^i_t(X^i_{t}, \mathbf U_{1:t-1}),  \quad i=1,2,
\end{equation}
where only actions $U^i_t \le X^i_t$ are feasible. Similar to Section \ref{Main_system_model}, we denote the  \emph{control strategy}  by $\bm g = (\bm g^1, \bm g^2)$.  The per unit cost $\ell(u^1_t,u^2_t)$ is defined to reflect the quality of transmission at time $t$ as follows:\vspace{-.0cm}
\begin{equation}
\ell(\mathbf x_t, \mathbf u_t)=\begin{dcases}
0 &  u^1_t=0,u^2_t=0\\
\ell^1\leq 0 &  u^1_t=1,u^2_t=0\\
\ell^2 \leq 0&  u^1_t=0,u^2_t=1\\
\ell^3 &  u^1_t=1,u^2_t=1\\
\end{dcases}\vspace{-.2cm}
\end{equation}
where $|\ell^j|\leq L,j=1,2,3$. The performance of strategy $\bm g$ is measured by \vspace{-.3cm}
\begin{equation}
   J(\bm g) = \mathbb{E}^{\bm g} \Big[
    \sum_{t=1}^\infty  \beta^{t-1} \ell(\mathbf X_t, \mathbf U_t) \Big].\vspace{-.2cm}
    \label{eq:discounted-reward}
\end{equation}
where $\beta \in (0,1)$. The case of symmetric arrivals ($p^1 = p^2)$ was considered
in~\cite{HluchyjGallager:1981, OoiWornell:1996}. In recent years, the above
model has been used as a benchmark for decentralized stochastic control problems~\cite{MNT:tractable-allerton,SeukenZilberstein:2007, DibangoyeMouaddibCahibdraa:2008}. We are interested in the following problem.
\begin{problem}\label{Problem2}
Given any  $\epsilon>0$, without knowing the arrival probabilities $p^1$ and $p^2$, and cost functions $\ell^1,\ell^2,\ell^3$, develop a decentralized Q-learning algorithm for both users such that users consistently learn an $\epsilon$-optimal strategy  $\bm g^\ast$.
\end{problem}

\subsection{Decentralized Q-learning Algorithm}
In this section, we follow the proposed two-step approach to develop a finite-state RL algorithm.
\subsubsection{An Equivalent Centralized POMDP}
In this step, we follow \cite{Mahajan2010MABC} and obtain the equivalent centralized  POMDP  for the completely known model as described in Section \ref{CIA_general}.

The common information shared between users is~$C_t=\mathbf{U}_{1:t-1}$. Define~ $Z_t=C_{t+1} \backslash C_t=\mathbf U_t$. At time~$t$, the coordinator observes~ $C_t=Z_{1:t-1}$ and prescribes~$\gamma^i_t \colon X^i_t \mapsto U^i_t$ that tell each agent how to use their local information to generate the control action.  For this specific model, the prescription $\gamma^i$ is completely
specified by~$A^i_t \coloneqq \gamma^i_t(1)$ (since~$\gamma^i_t(0)$ is always~$0$).
Hence, \vspace{-.1cm}
\begin{equation}
  U^i_t = \gamma^i_t(X^i_t) = A^i_t \cdot X^i_t\vspace{-.1cm}
  \label{eq:prescription}
\end{equation}
Therefore, we may equivalently assume that the coordinator generates actions
$\mathbf A_t = (A^1_t, A^2_t)$. The agents are passive and generate
actions $(U^1_t, U^2_t)$ according to~\eqref{eq:prescription}. Hence, at time $t$, the coordinator  prescribes action $\mathbf A_t \in \{0,1\}^2$ and observes $Z_t=\mathbf{U}_t\in \{0,1\}^2$.

Following \cite{Mahajan2010MABC}, define $\bm \Pi_t = (\Pi^1_t, \Pi^2_t)$, $\Pi^i_t = \mathds P(X^i_t = 1 \mid
\mathbf U_{1:t-1},\mathbf{A}_{1:t-1})$, as information state for the coordinated system with initial state $\bm \Pi_1=(p^1,p^2)$. It is shown in \cite{Mahajan2010MABC}:
\vspace{.1cm}
\hspace{-.1cm}
1) The information state $\bm \Pi_t$ evolves according to
 \begin{equation}
 \bm \Pi_{t+1}=\phi(\bm \Pi_t, \mathbf A_t,\mathbf U_t)
 \end{equation}
 where 
\begin{equation}
\hspace*{-.2cm}\phi(\bm \Pi_t, \hspace{-.05cm}\mathbf A_t,\hspace{-.1cm}\mathbf U_t\hspace{-.05cm})\hspace{-.1cm}=\hspace{-.1cm}
  \begin{cases}
\hspace{-.1cm}    (T_1\Pi^1_t,T_2\Pi^2_t) &  \mathbf{A}_t=(0,0)\\
\hspace{-.1cm}        (p^1,T_2\Pi^2_t)     &  \mathbf{A}_t=(1,0)\\
\hspace{-.1cm}    (T_1\Pi^1_t,p^2)     &  \mathbf{A}_t=(0,1)\\      
\hspace{-.1cm}    (1,1) &  \mathbf{A}_t\hspace{-.1cm}=\hspace{-.1cm}(1,1),  \mathbf{U}_t\hspace{-.1cm}=\hspace{-.1cm}(1,1)\\
\hspace{-.1cm}    (p^1,p^2)   &   \mathbf{A}_t\hspace{-.1cm}=\hspace{-.1cm}(1,1),  \mathbf{U}_t \hspace{-.1cm}\neq \hspace{-.1cm}(1,1)
  \end{cases}
\end{equation}
where $(p^1,p^2)$ are arrival rates and operator $T_i$ is given by
$
  T_i q = (1-p^i)(1-q), \quad i = 1,2.
$

\vspace*{.1cm}
\hspace{-.5cm}
2) The expected cost function is as follows:\vspace{-.1cm}
\begin{equation}
\hspace*{-.2cm} \hat{\ell}(\hspace{-.01cm}\bm \Pi_t\hspace{-.02cm}, \hspace{-.04cm}\mathbf A_t\hspace{-.02cm})\hspace{-.1cm}=\hspace{-.1cm}
  \begin{cases}
\hspace{-.08cm}    0, & \hspace{-.1cm} \mathbf A_t\hspace{-.1cm}=\hspace{-.1cm}(0,0)\\
\hspace{-.08cm}   \ell^1\Pi^1_t,           &  \hspace{-.1cm} \textbf{A}_t\hspace{-.1cm}=\hspace{-.1cm}(1,0)\\
 \hspace{-.08cm}   \ell^2\Pi^2_t,    & \hspace{-.1cm} \textbf{A}_t\hspace{-.1cm}=\hspace{-.1cm}(0,1)\\
 \hspace{-.08cm}  \ell^1 \hspace{-.05cm}\Pi^1_t\hspace{-.1cm}+\hspace{-.1cm}\ell^2\Pi^2_t\hspace{-.1cm}+\hspace{-.1cm}(\hspace{-.05cm}\ell^3\hspace{-.15cm}-\hspace{-.05cm} \ell^1\hspace{-.15cm}-\hspace{-.05cm}\ell^2)\Pi^1_t \hspace{-.05cm}\Pi^2_t  & \hspace{-.1cm} \textbf{A}_t\hspace{-.1cm}=\hspace{-.1cm}(1,1)\\
  \end{cases}
\end{equation} 
The action $(0,0)$ that corresponds to not transmitting is dominated by the actions $(1,0)$ or $(0,1)$. Therefore, with no loss of optimality, action $(0,0)$ is removed. In the sequel, we denote $\mathcal{A}:=\{(0,1),(1,0),(1,1)\}$ as the action space of the coordinator.

 We denote $\mathcal{R}$ as the reachable set of above centralized POMDP  that contains all the realizations of $\bm \pi_t$  generated by $\bm \pi_{t+1}=\phi(\bm \pi_t,\bm a, \bm u), \forall \bm a \in \mathcal{A}, \forall \mathbf{u} \in \{0,1\}^2,\forall t\in \mathbb{N}$,  with initial information state $\bm \pi_1 =(p^1,p^2)$. Thus, the reachable set  $\mathcal{R}$ is
given by \vspace{-.2cm}
\begin{multline}
\mathscr R \coloneqq \{ (1,1), (1,p^1), (p^2,1), (p^1,p^2) \} \\
  \cup \{ (p^1, T^n_2 p^2) : n \in \mathds N \}
  \cup \{ (T^n_1 p^1, p^2) : n \in \mathds N \},
\end{multline}
where
$
  T^n_i q = T_i ( T^{n-1}_i q).
$
According to Theorem \ref{Theorem1}, we have 
\begin{theorem}
Let $\bm \psi^\ast(\bm \pi)$ be any argmin of the right-hand side of the following dynamic program. For $\bm \pi \in \mathcal{R}$,
\begin{equation*} 
V(\bm \pi) \hspace{-.05cm}= \hspace{-.05cm}\min_{\bm a}(\hat{\ell}(\bm \pi \hspace{-.05cm}, \hspace{-.05cm}\bm a)+\beta \Exp{V(\phi(\bm \pi,\bm a,\mathbf U_{t}))|\bm \Pi_t \hspace{-.05cm}= \hspace{-.05cm}\bm \pi,\bm A_t \hspace{-.05cm}= \hspace{-.05cm}\bm a})
\end{equation*}
where $\bm a \in \mathcal{A}$. The stationary strategy $\bm g^\ast=(g^{1,\ast},g^{2,\ast})$ is optimal such that 
\begin{equation*}
g^{i,\ast}( \bm \pi, x)=\bm \psi^{i,\ast}(\bm \pi)\cdot x, \quad \forall \bm \pi \in \mathcal{R}, x \in \{0,1\}, i=1,2
\end{equation*}
where $\bm \psi^{i,\ast}$ denotes $i$th term of $\bm \psi^\ast$.
\end{theorem}

\begin{figure}[t!]
\centering
\vspace{-.4cm}
\hspace{-4cm}
\scalebox{.6}{
\leftline{\includegraphics[trim=0cm 5cm 0cm 11cm,clip, width=15cm]{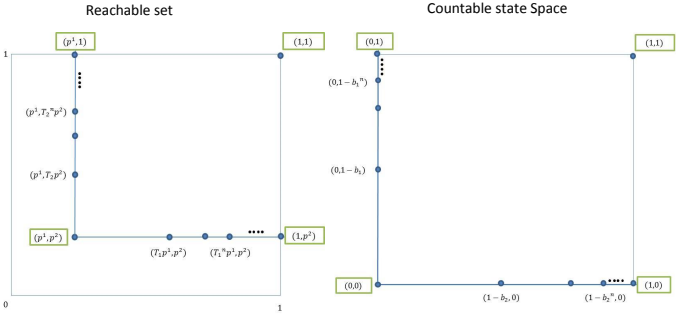}} }
\vspace{-.5cm}
\caption{It shows the reachable set $\mathcal{R}$ and the countable state space $\mathcal{S}$.}
\vspace*{-.5cm}
\end{figure}

\subsubsection{Q-learning algorithm for the  POMDP}

Let $b_1,b_2$ be any arbitrary number in $(0,1)$ and $B:\mathcal{R} \mapsto {\mathbb{Z}^+}^2$ be a bijective function that maps each state of $\mathscr R$  to a point in ${\mathbb{Z}^+}^2$ as follows:
\[\left(0,1-b_1^{n}\right)\hspace{-.1cm}=\hspace{-.1cm}B^{}(p^1,T_2^np^2),
\left(1-b_2^{n},0\right)\hspace{-.1cm}=\hspace{-.1cm}B^{}(T_1^np^1,p^2), n\in \mathds{N}
\]
\[\left(0\hspace{-.02cm},\hspace{-.05cm}1\right)\hspace{-.1cm}=\hspace{-.1cm}B^{}(p^1\hspace{-.05cm},\hspace{-.05cm}1)\hspace{-.02cm},\hspace{-.02cm} (1,0)\hspace{-.1cm}=\hspace{-.1cm}B^{}(1\hspace{-.05cm},\hspace{-.05cm}p^2), \left(1\hspace{-.05cm},\hspace{-.05cm}1\right)\hspace{-.1cm}=\hspace{-.1cm}B^{}(1\hspace{-.05cm},\hspace{-.05cm}1)\hspace{-.03cm}, \hspace{-.03cm} (0\hspace{-.05cm},\hspace{-.05cm}0)\hspace{-.1cm}=\hspace{-.1cm}B^{}(p^1\hspace{-.05cm},\hspace{-.05cm}p^2\hspace{-.03cm})
\] 
where $\lim_{n \rightarrow \infty} B(p^1,T_2^np^2)=B(p^1,1)$ and $\lim_{n \rightarrow \infty} B(T_1^np^1,p^2)=B(1,p^2)$.  Define a countable-state MDP $\Delta$ with state space $\mathcal{S}$, action space $\mathcal{A}$, dynamics $\tilde{f}$, and cost function $\tilde{\ell}$ as follows:

\textbf{(F1)}  Let $\mathcal{S}=\{S_k\}_{k=1}^\infty$ be the state space, where $S_{1}=\{(0,0)\}$ and   $S_{k}=\{(0,0),(0,1),(1,0),(1,1),(0,1-b_1^{i}),(1-b_2^{i},0)\}_{i=1}^{k-1}$, $k\geq 2$. The action space is $\mathcal{A}=\{(0,1),(1,0),(1,1)\}$. The initial state $S_1=(0,0)$.  The state $S_t \in \mathcal{S}_k$, $k \leq t$, evolves as follows: for $\mathbf A_t \in \mathcal{A}, \mathbf U_t \in \{0,1\}^2$,
\begin{equation}
S_{t+1}=\tilde{f}(S_t,\mathbf{A}_t,\mathbf U_t), \quad  S_{t+1} \in \mathcal{S}_{k+1}.
\end{equation}
For ease of exposition of dynamics $\tilde{f}$, we denote every state  $S_t \in \mathcal{S}_k$ in a format of $(1-b_2^{k_2},1-b_1^{k_1})$, where $k_1,k_2$ take value in the set of $\{0,1,\ldots,\infty\}$. Thus,
\begin{equation*}
\hspace{-.5cm}\tilde{f}(S_t,\mathbf{A}_t,\mathbf U_t)\hspace{-.1cm}=\hspace{-.1cm}
  \begin{cases}
\hspace{-.1cm}        (0,1-b_1^{(k_1+1)})     &  \mathbf{A}_t\hspace{-.1cm}=\hspace{-.1cm}(1,0)\\
\hspace{-.1cm}    (1-b_2^{(k_2+1)},0)    &  \mathbf{A}_t\hspace{-.1cm}=\hspace{-.1cm}(0,1)\\      
 \hspace{-.1cm}   (1,1) &  \mathbf{A}_t\hspace{-.1cm}=\hspace{-.1cm}(1,1),  \mathbf{U}_t\hspace{-.1cm}=\hspace{-.1cm}(1,1)\\
\hspace{-.1cm}    (0,0)   &   \mathbf{A}_t\hspace{-.1cm}=\hspace{-.1cm}(1,1),  \mathbf{U}_t \hspace{-.1cm}\neq \hspace{-.1cm}(1,1)
  \end{cases}
\end{equation*}
\hspace{-.05cm}
At time $t$,  there is a cost  given by
\begin{equation}\label{Countable_state_cost_MABC}
\tilde{\ell}(S_t,  \mathbf A_t)=\hat{\ell}(B^{-1}(S_t),\mathbf A_t).
\end{equation}

It is trivial to see that the tuple $\langle  
\{\mathcal{S}_k\}_{k=1}^\infty,B^{-1},\tilde{f}\rangle$ is an IER  because of the fact that \[\phi(\cdot,\mathbf a, \mathbf u)=B^{-1}\left(\tilde{f}\left(B(\cdot),\mathbf{a},\mathbf u) \right) \right), \quad \forall \mathbf a \in \mathcal{A}, \mathbf u \in \{0,1\}^2.\]

\vspace*{-.7cm}

\textbf{(F2)} State space $\mathcal{S}$, action space $\mathcal{A}$, and dynamics $\tilde{f}$ do not depend on the unknowns i.e. $(p^1,p^2,\ell^1,\ell^2,\ell^3)$.

The performance of a stationary strategy $ \tilde{\bm \psi}: \mathcal{S} \mapsto \mathcal{A}$ is quantified by  \eqref{Delta_total_cost}.
According to Lemma \ref{Delta_lemma}, we can restrict attention in solving MDP $\Delta$ instead of the POMDP without loss of optimality.  Let $\Delta_N$ be a finite-state MDP with state space  $\mathcal{S}_N$ and action space $\mathcal{A}$. The initial state  $S_1=(0,0)$. At time $t$,  state $S_t \in \mathcal{S}_N$ evolves as follows: for any $\mathbf A_t \hspace{-.05cm} \in \hspace{-.05cm} \mathcal{A}, \mathbf U_t \hspace{-.05cm} \in \hspace{-.05cm} \{0,1\}^2$,
\begin{equation}\label{Dynamics_delta_N_MABC}
S_{t+1}\hspace{-.05cm}=\hspace{-.05cm}\begin{dcases}
\hspace{-.05cm}\tilde{f}(S_t,\mathbf A_t,\mathbf U_t) & \tilde{f}(S_t,\mathbf A_t,\mathbf U_t)\hspace{-.05cm}  \in \hspace{-.05cm}\mathcal{S}_N\\
\hspace{-.05cm} (1-b_2,0) & \tilde{f}(S_t,\mathbf A_t,\mathbf U_t) \hspace{-.05cm} \in \hspace{-.05cm} \mathcal{S}_{N+1}\hspace{-.05cm}\backslash\hspace{-.01cm} \mathcal{S}_N 
\end{dcases}
\end{equation}
In \eqref{Dynamics_delta_N_MABC}, whenever  state $s_t$ steps out of $S_N$, the users take a sequence of actions as follows: At first, user 1 transmits and user 2 does not transmit, then user 2 transmits and user 1 does not transmit. This sequence of actions takes the system to state $(1-b_2,0) \in S_N, N \geq 2$. Now, we use standard Q-learning algorithm as the generic RL algorithm $\mathcal{T}$ to learn the optimal strategy of $\Delta_N$.

\vspace{-0.0cm}%
\alglanguage{pseudocode}
\begin{algorithm}[t!]
\small
\caption{\hspace{0.3cm}Decentralized Q-learning Algorithm}
\label{Q-learning_algorithm_MABC}
\begin{algorithmic}[1]

\State  Given $\epsilon>0$, choose a sufficiently large $N \in \mathbb{N}$ such that $\frac{2\beta^N}{1-\beta}L \leq \epsilon$. Then, construct state space $\mathcal{S}_N$, action space $\mathcal{A}$, and dynamics $\tilde{f}$. Let $s_1=(0,0)$. Initialize  Q-functions with  zero and  step-sizes  $\alpha$ with one i.e $Q(s,\mathbf{a})=0,\alpha(s,\mathbf{a})=1,  \forall s \in \mathcal{S}_N, \forall \mathbf{a}\in \mathcal{A}$. 

\State    At iteration $k \in \mathbb{N}$, users uniformly  pick a random action  $\ \mathbf a_k \in \mathcal{A}$ at state $s_k \in \mathcal{S}_N$ by means of a common shared  random number generator. Then, user $i\in \{1,2\}$ takes action $u^i_k$ according to the chosen  $a^i_k$ and local information $x^i_k \in \{0,1\}$ as follows:\vspace{-.2cm}
\begin{equation*}
u^i_k=a^i_k\cdot x_k^i. \quad  i=1,2.\vspace{-.1cm}
\end{equation*}

\State  Based on the taken actions, the system incurs a cost $\ell_k$ and   generates  $(x^1_{k+1},x^2_{k+1},\mathbf u_k=(u^1_k,u^2_k))$. Since $\mathbf{u}_k$ is observable to both users, they  consistently compute the next state\vspace{-.1cm}
\begin{equation*}
 s_{k+1}=\tilde{f}(s_k,\mathbf  a_k,\mathbf u_k).\vspace{-.1cm}
\end{equation*}
 If $s_{k+1} \notin  \mathcal{S}_N$,  user 1 transmits first and then, user 2 transmits, and  the state of system will be transmitted to $s_{k+1}=(1-b_2,0)$; otherwise, the system proceeds from $s_{k+1} \in \mathcal{S}_N$.
 
 \State Users update the corresponding Q-function associated with the pair $(s_k,\mathbf a_k)$ as follows:
 \vspace{-.2cm}
\hspace{-.1cm}
\begin{equation*}
\hspace{-.2cm}Q(\hspace{-.05cm}s_k\hspace{-.05cm},\hspace{-.05cm}\mathbf{a}_k\hspace{-.05cm})\hspace{-.05cm}\leftarrow \hspace{-.05cm}(\hspace{-.05cm}1-\alpha(\hspace{-.05cm}s_k\hspace{-.05cm},\hspace{-.05cm}\mathbf{a}_k\hspace{-.05cm})\hspace{-.05cm})Q(\hspace{-.05cm}s_k \hspace{-.05cm}, \hspace{-.05cm}\mathbf{a}_k\hspace{-.05cm})+
\alpha(\hspace{-.05cm}s_k\hspace{-.05cm},\hspace{-.05cm}\mathbf{a}_k\hspace{-.05cm})\hspace{-.1cm}\left(\hspace{-.1cm}\ell_k\hspace{-.1cm}+\hspace{-.1cm}\beta\hspace{-.05cm} \min_{v \in \mathcal{A}}Q(\hspace{-.05cm}s_{k+1}\hspace{-.05cm},\hspace{-.05cm}v\hspace{-.05cm}) \hspace{-.1cm}\right)\vspace{-.1cm}
\end{equation*} 
\vspace{-.1cm}
Also, the corresponding step-size are updated: \vspace{-.1cm}
 \[\frac{1}{\alpha(s_k,\mathbf{a}_k)} \leftarrow \frac{1}{\alpha(s_k,\mathbf{a}_k)} +1.\]
  \State \noindent $k \leftarrow k+1$, and go to step 2 until termination.
\Statex
\end{algorithmic}
  \vspace{-0.4cm}%
\end{algorithm}
\vspace{-.0cm}%

According to~\cite[Theorem 3]{Tsitsiklis1994asychronous}, Q-functions in Algorithm \ref{Q-learning_algorithm_MABC} will  converge to a $Q^\ast$ with probability one\footnote{In this example, every pair of (state,action) will be visited infinitely often  by uniformly randomly picked actions.}.  Let $Q^\ast$ be the resultant limit. Then,  the optimal strategy $\tilde{\bm \psi_N^\ast}$  is  as follows:
\begin{equation}\label{greedy}
\tilde{\bm \psi_N^\ast}= \argmin_{\mathbf a \in \mathcal{A}}(Q^\ast(\cdot,\mathbf a)).
\end{equation}
The strategy $\bm g^\ast_N$ is $\epsilon_N$-optimal where
\begin{equation*}
g^{i,\ast}_N(s)(x):=\tilde{\bm \psi_N^{i,\ast}}(s)\cdot x, \quad \forall s \in \mathcal{S}_N, x \in \{0,1\}, i=1,2
\end{equation*}
where $\tilde{\bm \psi_N^{i,\ast}}$ denotes $i$th term of $\tilde{\bm \psi_N^{\ast}}$.

\subsubsection{Numerical Results}
In this section, we provide a numerical simulation that shows the strategy learned by  decentralized Q-learning  Algorithm \ref{Q-learning_algorithm_MABC} converges to an optimal strategy when the arrival probabilities are $(p^1,p^2)=(0.3,0.6)$ and the cost functions are $\ell^1=\ell^2=-1, \ell^3=0$. 

\begin{figure}[t!]
\centering

\vspace{-0cm}
\hspace{-0cm}
\scalebox{1.05}{
\leftline{\includegraphics[trim=1.5cm 0cm 0cm 0cm,clip, width=8cm]{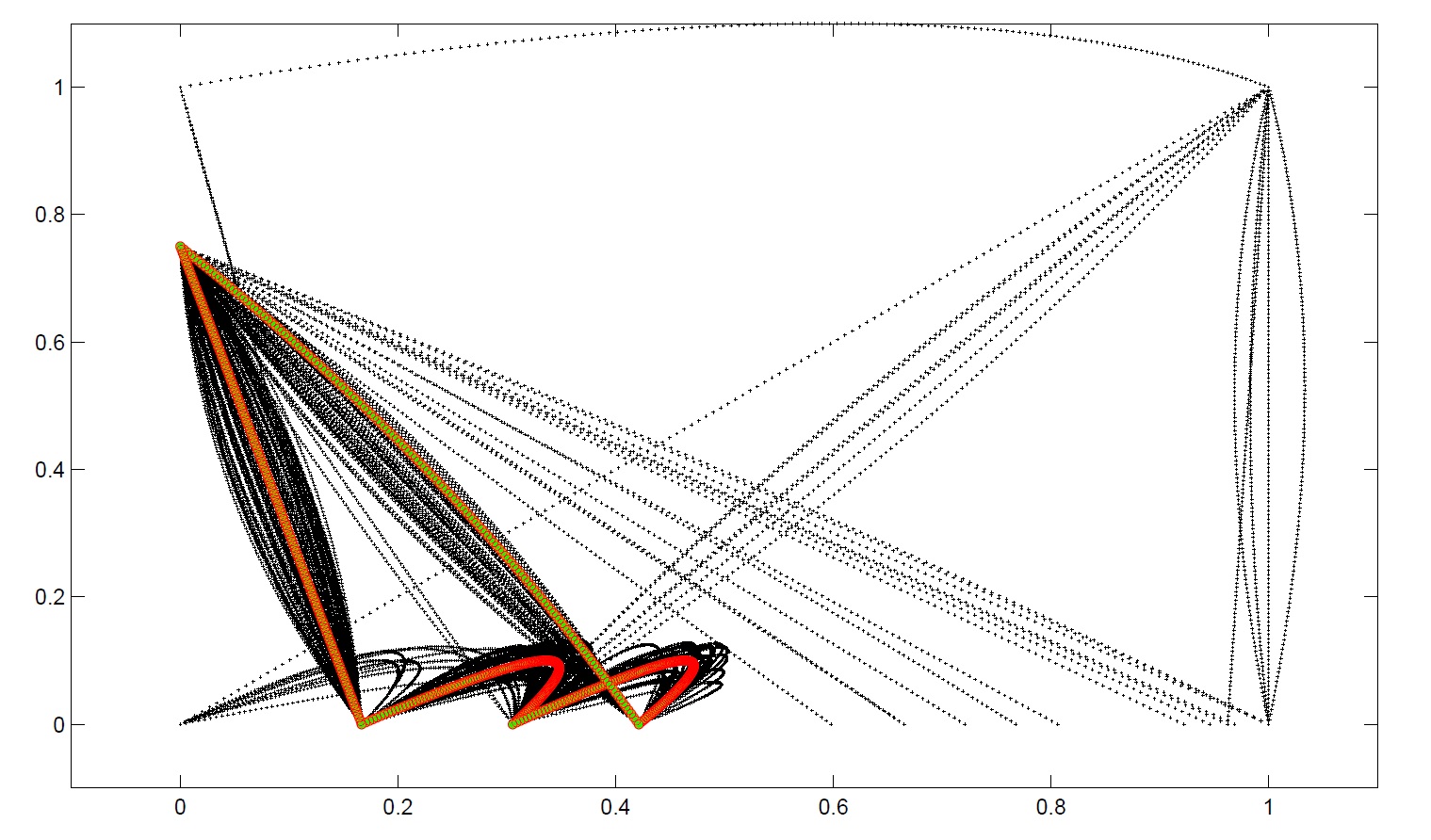}}}
\vspace{-.7cm}
\caption{This  figure displays the learning procedure of optimal strategy in a few snapshots. It is seen that the state
of the system is eventually trapped in the optimal recurrent class. The learning procedure is plotted in black and the
optimal recurrent class is plotted in red. In this simulation, we use the following numerical values: $b_1=0.25,b_2=0.83,N=20, \beta=0.99, p^1=0.3,p^2=0.6, \ell^1=\ell^2=-1,\ell^3=0.$}
\label{Si}
\vspace{-.6cm}
\end{figure}

Suppose users have no packets at the beginning. Users wait one time step  to receive packets (i.e. user 1 receives a packet with $p^1=0.3$ probability and user 2 receives a packet with $p^2=0.6$ probability). At $t=1$, action $(0,1)$ is optimal i.e. the user $2$ transmits and user $1$ does not transmit. At $t\geq 2$, state $s_t$ enters a recurrent class  under the optimal strategy, and stays there forever. 
The recurrent class includes four states:  $(0,1-b_1^1)$, $(1-b_2^1,0)$, $(1-b_2^2,0)$, and $(1-b_2^3,0)$. One immediate result is that for any $N \geq 4$, state $S_t$ will never step out of $\mathcal{S}_N$ under the optimal strategy which implies $\tau_N=\infty$ and hence $\epsilon_N$ in Theorem \ref{Theorem3} is zero (i.e. optimal strategy). For $t\geq 2$, the optimal strategy is a sequence of the following actions $(0,1),(0,1),(1,0),(0,1)$. Thus, it means that user 2 should transmit $3$ times more than user 1 to  minimize the number of collisions (maximize the number of successful transmission). Figure \ref{Si} displays  a few snapshots of state $s_t$ governed by the strategy under the learning procedure where the learned strategy will eventually  take state $s_t$ to the optimal recurrent class. 

\vspace{-.1cm}

%

\section{Conclusion}

In this paper, we proposed a novel approach to develop a finite-state RL algorithm, for a large class of decentralized control systems with partial history sharing information structure, that guarantees $\epsilon-$team-optimal solution. We presented our approach  in two  steps. In the first step, we used the common information approach  to obtain an equivalent centralized POMDP of the decentralized control problem. However, the resultant POMDP can not be used directly because it requires the complete knowledge of the model while the agents only know the model incompletely. Thus,  in the second step, we introduced a new methodology to develop a RL algorithm for the obtained centralized POMDP. In particular, to remove the dependency of the complete knowledge, we introduced  Incrementally Expanding Representation (IER) and  based on that, we constructed a finite-state RL algorithm. In addition, we illustrated our approach by developing a decentralized Q-learning algorithm for two-user Multi Access Broadcast Channel (MABC), a benchmark example for decentralized control systems. The numerical simulations verify that the learned strategy converges to an optimal strategy.

\vspace{-.2cm}
\bibliographystyle{IEEEtran}

\end{document}